\documentclass[11pt]{article}
\usepackage[margin=1in]{geometry}
\usepackage{amsmath,amssymb,amsthm,mathtools}
\usepackage[hidelinks]{hyperref}
\usepackage{xcolor}

\newtheorem{theorem}{Theorem}[section]
\newtheorem{proposition}{Proposition}[section]
\newtheorem{lemma}{Lemma}[section]
\newtheorem{corollary}{Corollary}[section]
\newtheorem{theoremA}{Theorem}

\newtheorem{lemmaA}{Lemma}

\theoremstyle{definition}
\newtheorem{definition}{Definition}[section]
\theoremstyle{remark}
\newtheorem{remark}{Remark}[section]

\newcommand{\Sd}{\mathbb{S}^d}
\newcommand{\Pt}{\mathcal{P}_t}
\newcommand{\dd}{\,\mathrm{d}}
\newcommand{\dist}{\operatorname{dist}}
\newcommand{\Lam}{\Lambda}

\title{Nested Spherical Designs of Optimal Order\thanks{Department of Mathematics, City University of Hong Kong, Hong Kong. Emails: alirazajatoi1@gmail.com (A. R. Jatoi); xzhuang7@cityu.edu.hk (X. Zhuang). } }
\author{Ali Raza Jatoi
\and Xiaosheng Zhuang}

\date{}

\begin{document}
\maketitle
\begingroup
\renewcommand{\thefootnote}{}
\footnotemark
\endgroup
\begin{abstract}

Zheng and Zhuang \cite{ZZ} showed that a spherical $t_1$-design $Y$ with $|Y|=M\sim t_1^d$ points on the sphere $\Sd$ can be extended to a spherical $t$-design by adding points, with a total of order $t^{2d+1}$, and asked whether the optimal order $t^d$ is attainable. We keep their construction and improve the estimate of the pairing of the given points against polynomials of unit $L^1$-gradient norm. We prove that, if the given spherical $t_1$-design is upper regular at scale $t_1$ (in particular if it is well-separated; such a design exists by Bondarenko, Radchenko and Viazovska \cite{BRV15}), then for $d\ge2$, every $N\ge C(d,K)\,t^d$ suffices to admit a point set $X$ with $|X|=N$, for all $t_1<t$, such that the multiset union $X\uplus Y$ is a spherical $t$-design; in particular $|X\uplus Y|\le C(d,K)\,t^d$. The proof combines the design property, which removes the low frequencies exactly, with an $L^1$ Jackson-type estimate and a sampling inequality applied on dyadic frequency blocks. Upper regularity turns out to be automatic for nested spherical designs of optimal order; that is, every spherical $t_1$-design with $M\le C_M t_1^d$ points extends to a spherical $t$-design with $N+M\le C(d)\,C_M\,t^d$ (with $C_M\ge1$), uniformly in $t_1<t$, without any separation hypothesis. Some control of local multiplicity is necessary for an extension of optimal order. We also show that the estimate of the pairing functional cannot be improved within this method; that is, for every point set of at most $c\,t^d$ points it is at least of order $t^{d-1}$, and for upper regular designs it is of exactly this order. Finally, we treat finite chains of nested spherical designs, and on the circle we remove a logarithmic factor from the general estimate.
\end{abstract}


\section{Introduction}\label{sec:intro}

Let $\Sd=\{x\in\mathbb{R}^{d+1}:|x|=1\}$ denote the normalized unit sphere in $\mathbb{R}^{d+1}$ with the normalized Lebesgue surface measure $\mu$ such that $\mu(\Sd)=1$, and let $\Pi_t:=\Pi_t(\Sd)$ denote the space of all polynomials of degree
at most $t$ restricted on $\Sd$. A finite multiset $X=\{x_1,\dots,x_N\} \subset \Sd$ is a \emph{spherical $t$-design} if the equal-weight rule on $X$ integrates every polynomial of degree at most $t$
exactly:
\[
  \frac1N\sum_{i=1}^NQ(x_i)=\int_{\Sd}Q\dd\mu\qquad\text{for all }Q\in\Pi_t(\Sd).
\]
The notion was introduced by Delsarte, Goethals and Seidel \cite{DGS77}, building on earlier work on invariant cubature, and has since connected approximation theory, algebraic combinatorics, coding theory, experimental design and, more recently, learning on spheres; see \cite{AnZhuang} for a recent overview. Because the weights are equal, a spherical design is a particularly economical quadrature rule. It is stable and requires no weight storage, and a spherical $2L$-design can serve as the quadrature rule in hyperinterpolation of degree $L$ \cite{Sloan95}, which in general uses positive but not necessarily equal weights.

Let $N(d,t)$ be the minimal cardinality of a spherical $t$-design on $\Sd$. The bound of Delsarte--Goethals--Seidel \cite{DGS77} gives $N(d,t)\ge c\,t^d$, because the space $\Pi_t$ has dimension of order $t^d$; designs meeting the bound (\emph{tight designs}) exist for every $t$ when $d=1$ (the regular $(t+1)$-gon), but for $d\ge2$ only if $t\in\{1,2,3,4,5,7,11\}$, and then only in very few cases \cite{BB09,AnZhuang}. Existence of spherical $t$-designs for all $d$ and $t$ was proved by Seymour and Zaslavsky \cite{SZ84}, non-constructively and without a bound on the number of points. Subsequent upper bounds were polynomial in $t$ but far from the lower bound: $t^{Cd^4}$ (Wagner--Volkmann \cite{WV91}), $t^{Cd^3}$ (Bajnok \cite{Baj92}), $t^{(d^2+d)/2}$ (Korevaar--Meyers \cite{KM93}), and $t^{2d(d+1)/(d+2)}$ (Bondarenko--Viazovska \cite{BV10}, via Brouwer's fixed point theorem). The long-standing conjecture $N(d,t)\le C_d\,t^d$ was settled by Bondarenko, Radchenko and Viazovska \cite{BRV13}, who combined Marcinkiewicz--Zygmund inequalities \cite{MNW01} with Brouwer degree theory on the space of polynomials; in a second paper they showed that the designs can in addition be taken to be \emph{well-separated} \cite{BRV15}. Thus $t^d$ is the optimal order, attained by designs whose points are well-separated at the scale $N^{-1/d}$.

On the numerical side, spherical designs on $\mathbb S^2$ with $N$ close to the conjectured optimal order $t^2/2$ have been computed by optimization methods, for $t$ up to $1000$ in \cite{GrafPotts}, $t$ up to $3200$ in \cite{XiaoZhuang2025}, and in the tables of \cite{HardinSloane,Womersley}; \cite{Womersley} also explores $d=3$. This makes the question of \emph{reusing} an existing spherical design practically relevant.

In multilevel or adaptive computations, a quadrature rule is refined repeatedly for a fixed integrand. If the node sets are nested, earlier function evaluations are reused, and only the new nodes require evaluation. Zheng and Zhuang \cite{ZZ} introduced the corresponding notion for spherical designs. Given a spherical $t_1$-design $Y$ with $M$ points and $t>t_1$, find a set $X$ of $N$ points such that the (multiset) union $X\uplus Y$ is a spherical $t$-design; then $Y\subset X\uplus Y$ is a \emph{nested} pair of spherical designs. They asked three questions: (Q1) Do nested spherical designs exist? (Q2) How many points are needed to extend a spherical $t_1$-design to a spherical $t$-design? (Q3) What is the optimal order?

Zheng and Zhuang answered (Q1) affirmatively, in fact, for an arbitrary prescribed point set $Y\subset \Sd$, by adapting the Brouwer-degree argument of \cite{BRV13}. For (Q2), when $Y$ is a spherical $t_1$-design of optimal order $M\sim t_1^d$, they showed that the total number of points $N+M$ needed is of order $t^{2d+1}$ \cite{ZZ}, which exceeds the optimal order $t^d$ by a factor $t^{d+1}$. They describe this estimate as pessimistic and formulate a conjecture about the constants in the optimal orders for spherical $t_1$- and $t$-designs. Their concluding remarks also identify two directions: to find out whether properties such as being well-separated improve the estimate, and the role of a quantity measuring how far a given spherical design $Y$ is from being a good starting point.

That quantity appears explicitly in our analysis. Let $\Pt:=\Pt(\Sd)=\{P\in \Pi_t: \int_{\Sd} P\dd\mu=0\}$ be the Hilbert space of mean-zero polynomials on $\Sd$ and $\Omega:=\{P\in\Pt: \int|\nabla P|\dd\mu<1\}$ so that $\partial\Omega=\{P:\int|\nabla P|\dd\mu=1\}$. Let $G_y\in\Pt$ denote the reproducing kernel at $y$. For a finite multiset $Y\subset\Sd$, define 
\begin{equation}\label{eq:Lambda}
  \Lam(Y,t):=\sup_{P\in\partial\Omega}\Bigl(-\sum_{y\in Y}P(y)\Bigr)_+
  =\sup_{P\in\partial\Omega}\Bigl(-\Bigl\langle\sum_{y\in Y}G_y,\,P\Bigr\rangle\Bigr)_+ .
\end{equation}
The degree-theoretic argument of \cite{BRV13,ZZ} (see also Proposition~\ref{prop:master}) shows that $N\ge C_{1,d}t^d$ new points suffice as soon as $N\,c_d/t\ge\Lam(Y,t)$, where $c_d$ is the constant of \cite[Theorem~1]{ZZ}. The task is therefore to bound
$\Lam(Y,t)$. The estimate of \cite{ZZ} amounts to the bound $\Lam(Y,t)\lesssim M\,t^d$ and does not use the design property of $Y$. It is not sharp even when all points of $Y$ coincide. For $d\ge2$, the Sobolev and Nikolskii inequalities give
$\|P\|_\infty\le C\,t^{d-1}\int|\nabla P|\dd\mu$ on $\Pt$, so $\Lam(Y,t)\le C\,M\,t^{d-1}$ for every multiset $Y$ of $M$ points (see Remark~\ref{rem:SN}). The argument of Section~\ref{subsec:sharpness} gives a lower bound of the same order when all $M$ points coincide (see Remark~\ref{rem:meaning}).

Our key result (see Proposition~\ref{prop:dyadic}) improves the estimate of  $\Lambda(Y,t)$ by exploiting the design property of $Y$. We use a concept of upper regularity that is more general than well-separatedness, defined as follows. 

\begin{definition}\label{def:upper-regularity}
A finite multiset $Y\subset \Sd$ is \emph{upper regular at scale $t_1$ with a constant $K\ge1$} if every spherical cap of radius $r>0$ contains at most $K(1+rt_1)^d$ points of $Y$. 
\end{definition}
Every (well-separated) spherical $t_1$-design with a minimal separation distance $\ge c/t_1$ has this upper regularity property (see Lemma~\ref{lem:sep}), and such well-separated spherical designs with $M\sim t_1^d$ exist by \cite{BRV15}. With the upper regularity property, we show that $\Lambda(Y,t)\le  A_*t^{d-1}$ for some constant $A_*=A_*(d,K)$ depending only on $d$ and $K$, and it leads to the following main theorem for nested spherical designs of optimal order.

\begin{theorem}[Upper regular extension]\label{thm:extension}
Let $d\ge2$, $1\le t_1<t$, and let $Y=\{y_1,\ldots,y_M\}\subset \Sd$ be a spherical $t_1$-design that is upper regular at scale $t_1$ with a constant $K\ge1$. Then,
\[
   \Lam(Y,t)\le A_*\,t^{d-1},
\]
for some constant $A_*=A_*(d,K)$, and every $N\ge C(d,K)\,t^d$ admits a point set $X\subset\Sd$ with $|X|=N$ such that $X\uplus Y$ is a spherical $t$-design. If $|Y|=M\le C_M t_1^d$ with $C_M\ge1$, then
$N+M\le C(d,K,C_M)\,t^d$, which is of optimal order, uniformly in $t_1<t$.
\end{theorem}

The proof (see Section~\ref{subsec:extension}) uses the idea of dyadic refinements. For $d=1$ the same argument gives an additional factor $1+\log(t/t_1)$ for $\Lambda(Y,t)$, which Theorem~\ref{thm:circle} removes for spherical designs with $M\lesssim t_1$.   Furthermore, we show that the estimate $\Lambda(Y,t)\le A_* t^{d-1}$ is sharp, and the criterion of Proposition~\ref{prop:master} cannot certify fewer than order $t^d$ new points. 

Iterating Theorem~\ref{thm:extension} produces a chain $Y_0\subset Y_1\subset\cdots\subset Y_{L-1}$ of nested spherical designs, where $Y_k$ is a spherical $t_{k+1}$-design. Compared to \cite[Proposition 1]{ZZ} on the special case of spherical designs with optimal order, here, in our paper, we do not need the condition on the ratios $t_{k+1}/t_k$. The proof uses the fact that the points constructed in \cite[Appendix~B]{ZZ} are displaced from the cells of an area-regular partition by at most $r_d/(3t)$, so that upper regularity passes from $Y_{k-1}$ to $Y_k$ with a controlled constant. The constant grows geometrically in the number of levels $L$ (see Theorem~\ref{thm:chain}).




Upper regularity is in fact automatic for spherical designs of optimal order (Proposition~\ref{prop:auto}), so the hypothesis of Theorem~\ref{thm:extension} can be dropped.

\begin{theorem}[General extension]\label{thm:arbitrary}
Let $d\ge2$, $1\le t_1<t$, and let $Y=\{y_1,\ldots, y_M\}\subset \Sd$ be a spherical $t_1$-design with $M\le C_Mt_1^d$, where $C_M\ge1$. Then
\[
\Lam(Y,t)\le A_*(d,C_M)\,t^{d-1},
\]
and every $N\ge C(d)\,C_M\,t^d$ admits a point set  $X\subset \Sd$ with $|X|=N$
such that $X\uplus  Y$ is a spherical $t$-design. In particular, $N+M\le C(d)\,C_M\,t^d$,
uniformly in $t_1<t$.
\end{theorem}
Here $A_*(d,C_M)=\alpha_d\,K(d,C_M)$, where $\alpha_d$ is the constant of Section~\ref{subsec:chains} (with $A_*(d,K)=\alpha_dK$) and $K(d,C_M)$, linear in $C_M$, is the constant of Proposition~\ref{prop:auto}.

In summary, the comparison with \cite{ZZ} is as follows: Theorem~\ref{thm:arbitrary} improves their total from $t^{2d+1}$ to $t^d$ for every starting spherical $t_1$-design of optimal order, with no further hypothesis. It does not address their conjecture on constants. Zheng and Zhuang also observe \cite[Prop.~1]{ZZ} that for $t=mt_1$ with $m>1$ rational, the question becomes trivial for a specially chosen starting design: a spherical $t$-design of optimal order is also a spherical $t_1$-design, and it can be extended to a larger spherical $t$-design by adding copies of a spherical $t$-design, with identical constants. This uses a spherical design adapted to $t$ and applies only to rational multiples. Theorem~\ref{thm:arbitrary} applies to every spherical $t_1$-design with $M\lesssim t_1^d$ points and to all $t_1<t$, without a separation or regularity assumption.

The remainder of the paper is organized as follows. Section~\ref{sec:prelim} collects the preliminaries. Section~\ref{sec:upperregular} then constructs nested spherical designs from point sets that are upper regular. It proves Theorem~\ref{thm:extension} (Section~\ref{subsec:extension}) and the sharpness of $\Lambda(Y,t)$. Finite chains of nested spherical designs can then be derived from Theorem~\ref{thm:extension} (Section~\ref{subsec:chains}). Section~\ref{sec:beyond} removes the upper regularity assumption for spherical designs: Section~\ref{subsec:arbitrary} proves the master estimate and the necessity statement, shows that spherical designs of optimal order are automatically upper regular, and deduces Theorem~\ref{thm:arbitrary}. Section~\ref{subsec:circle} treats the circle, including the numerical test for $d=1$. Section~\ref{sec:conclusion} concludes the paper and lists some open problems: the case of arbitrary spherical $t_1$-designs beyond optimal order, the conjecture on constants in \cite{ZZ}, and chains of unbounded length.

\section{Preliminaries}\label{sec:prelim}

In this section, we provide preliminaries on the degree-theoretic framework and the results of \cite{BRV13,ZZ} for the existence of spherical designs, the filtered kernels and the $L^1$ Jackson-type estimate, and the upper regular point sets with the sampling bound.

\subsection{Existing Key Results}\label{subsec:existing}

Let $\|\cdot\|_p$ denote the usual $L^p$-norm for $L^p(\Sd,\mu)$, the space of measurable functions for which the $p$-th power of the absolute value is Lebesgue integrable. Recall that $\Pt=\{P\in\Pi_t: \int_{\Sd}P\dd\mu=0\}$ is the Hilbert space of mean-zero polynomials, $\Omega:=\{P\in\Pt:\int_{\Sd}|\nabla P|\dd\mu<1\}$, and $\partial\Omega=\{P\in\Pt:\int_{\Sd}|\nabla P|\dd\mu=1\}$, where for $P,Q\in\Pt$, $\nabla P$ is its Riemannian gradient and we use $\langle P,Q\rangle:=\int_{\Sd} PQ\dd\mu$ for the inner product of $P$ and $Q$. Throughout the paper, multisets are used, with points counted with their multiplicities. The union $X\uplus Y$ denotes multiset union for $X$ and $Y$. By the Riesz representation theorem, for each point $x\in\Sd$, there exists a unique polynomial $G_x\in \Pt$ such that $\langle G_x, P\rangle=P(x)$ for all $P\in\Pt$. A multiset $\{x_1,\dots,x_N\}$ is a spherical $t$-design if and only if $\sum_iG_{x_i}=0$. We use the following two results, which follow the degree-theoretic argument of \cite{BRV13}.

\begin{theoremA}[Zheng--Zhuang {\cite[Thm.~1]{ZZ}}]\label{thm:ZZ1}
There are constants $C_{1,d},c_d>0$ such that for all $t$ and all $N\ge C_{1,d}t^d$, there is a continuous map $F:\Pt\to(\Sd)^N$ with $F(P)=(x_1(P),\dots,x_N(P))$ satisfying $\sum_{i=1}^NP(x_i(P))>N\,c_d\,t^{-1}$ for all $P\in\partial\Omega$.
\end{theoremA}

\begin{theoremA}[Brouwer degree theory, in the form used in \cite{BRV13}]\label{thm:brouwer}
If $f:\Pt\to\Pt$ is continuous with $\langle P,f(P)\rangle>0$ on $\partial\Omega$,
then $f$ has a zero in $\Omega=\{P\in\Pt:\int|\nabla P|\dd\mu<1\}$. 
\end{theoremA}

\begin{proposition}[Master criterion]\label{prop:master}
Let $Y=\{y_1,\dots,y_M\}\subset\Sd$ be any multiset and let $\Lam(Y,t)$ be as in \eqref{eq:Lambda}.
If $N\ge C_{1,d}t^d$ and $N\,c_d/t\ge\Lam(Y,t)$, then there is a multiset $X\subset\Sd$ with $|X|=N$ such that $X\uplus Y$ is a spherical $t$-design. Moreover, $X$ can be taken of the form $X=F(Q)=\{x_i(Q)\}_{i=1}^N$ for some $Q\in\Omega$, with $F$ as in Theorem~\ref{thm:ZZ1}.
\end{proposition}
\begin{proof}
Let $F$ be as in Theorem~\ref{thm:ZZ1} and $f(P)=\sum_iG_{x_i(P)}+\sum_jG_{y_j}$. Then $f$ is continuous and $\langle P,f(P)\rangle=\sum_iP(x_i(P))+\sum_jP(y_j)>Nc_dt^{-1}-\Lam(Y,t)\ge0$ on $\partial\Omega$. By Theorem~\ref{thm:brouwer}, there is $Q\in \Omega$ with $f(Q)=0$, i.e. $\sum_i G_{x_i(Q)}+\sum_jG_{y_j}=0$. Then, $X$ is given by $X=F(Q)=\{x_i(Q)\}_{i=1}^N$. 
\end{proof}

\subsection{Filtered Kernels and Jackson-Type Estimate}\label{subsec:kernels}

Write $P_{\ell,d}$ for the Gegenbauer (Legendre) polynomial normalized by $P_{\ell,d}(1)=1$ and $D(\ell,d)$ for the dimension of the space of spherical harmonics of degree $\ell$. For a smooth function $h:[0,\infty)\to\mathbb{R}$ and $n\ge1$, define the filtered kernel $\Phi^h_n$ and its associated projection operator $V_n^h$ as 
\[
  \Phi^h_n(\langle x,y\rangle)=\sum_{\ell\ge0}h(\ell/n)\,D(\ell,d)\,P_{\ell,d}(\langle x,y\rangle),
  \qquad
  V^h_nf(x)=\int_{\Sd}\Phi^h_n(\langle x,y\rangle)f(y)\dd\mu(y),
\]
respectively,  where $\langle x,y\rangle$ is the inner product of vectors $x,y\in\Sd$.  Note that $V^h_n$ modulates the degree-$\ell$ harmonic component of $f$ by $h(\ell/n)$. Because $\int_{\Sd} P_{\ell,d}(\langle x,y\rangle)\dd\mu(y)=\delta_{\ell,0}$, we have
\begin{equation}\label{eq:mass}
  \int_{\Sd}\Phi^h_n(\langle x,y\rangle)\dd\mu(y)=h(0).
\end{equation}
Throughout the paper, we fix two cut-offs for $h$. They are fixed once and for all and do not depend on $d$; all constants below may depend on $h_1,h_2$, and we do not record this dependence.
\begin{itemize}
  \item[(a)] $h=h_1$ is smooth with values in $[0,1]$, $h_1\equiv1$ on $[0,\tfrac12]$, and $\operatorname{supp}h_1\subset[0,1]$. Then $V_n:=V^{h_1}_n$ maps $L^p(\Sd,\mu)$ into polynomials of degree $<n$ (since $h_1(1)=0$), satisfies $V_nR=R$ for every polynomial $R$ of degree $\le n/2$, and $\int_{\Sd}\Phi^{h_1}_n(\langle \cdot, y\rangle)\dd\mu(y)=h_1(0)=1$.  
  
  \item[(b)] $h=h_2$ is smooth, $h_2\equiv1$ on $[0,1]$, and $\operatorname{supp}h_2\subset[0,2]$. Then $V^{h_2}_nR=R$ for every polynomial $R$ of degree $\le n$.
\end{itemize}
The plateau of $h_1$ on $[0,\frac12]$ is used in Section~\ref{subsec:extension} to make the telescoping sum terminate exactly. Both cut-offs are constant near $0$, so their even extensions to $\mathbb{R}$ are smooth, which is the form of hypothesis under which the localization estimate is usually stated.

\begin{lemmaA}[Localization; Narcowich--Petrushev--Ward {\cite{NPW06}}, see also {\cite[Ch.~2]{DX13}}]\label{lem:loc}
For $h\in\{h_1,h_2\}$ and every $k\ge1$ there is $C_k=C_k(d,h)$ with
\[
  |\Phi^h_n(\cos\theta)|\le C_k\,n^d(1+n\theta)^{-k},\qquad 0\le\theta\le\pi,\ n\ge1 .
\]
\end{lemmaA}

This is the standard localization estimate for smoothly cut-off zonal kernels. The constant $C_k$ depends on $d$, $h$ and $k$ but not on $n$. Then, we have the following Jackson-type estimate, which uses the invariance of the Liouville measure.

\begin{lemma}[$L^1$ Jackson-type estimate]\label{lem:jackson}
There is $C_J=C_J(d,h_1)$ such that for all $n\ge1$ and all $f\in C^1(\Sd)$, it holds that
\[
  \|f-V_nf\|_1\le\frac{C_J}{n}\int_{\Sd}|\nabla f|\dd\mu .
\]
\end{lemma}

\begin{proof}
Write $\Phi=\Phi^{h_1}_n$. By \eqref{eq:mass} with $h_1(0)=1$,
\[
  f(x)-V_nf(x)=\int_{\Sd}\Phi(\langle x,y\rangle)\bigl(f(x)-f(y)\bigr)\dd\mu(y).
\]
Fix $x$ and use geodesic polar coordinates: $y=\exp_x(\theta v)$ with $\theta\in[0,\pi]$ and $v$ a unit vector in $T_x\Sd$. Then $\dd\mu(y)=\omega_d\sin^{d-1}\theta\dd\theta\dd v$, where $\dd v$ is the normalized measure on the unit sphere of $T_x\Sd$ and $\omega_d>0$ depends only on $d$ (for $d=1$ the unit sphere of $T_x\mathbb S^1$ is two points). Let $\gamma_{x,v}$ be the unit-speed geodesic with $\gamma_{x,v}(0)=x$, $\gamma_{x,v}'(0)=v$. Then $y=\gamma_{x,v}(\theta)$ and
\[
  |f(x)-f(y)|\le\int_0^\theta|\nabla f(\gamma_{x,v}(s))|\dd s .
\]
Integrating over $x$ and using Tonelli,
\[
  \|f-V_nf\|_1\le \omega_d\int_0^\pi|\Phi(\cos\theta)|\sin^{d-1}\theta
  \Bigl[\int_0^\theta\!\!\iint|\nabla f(\gamma_{x,v}(s))|\dd v\dd\mu(x)\,\dd s\Bigr]\dd\theta .
\]
The measure $\dd v\dd\mu(x)$ is the normalized Liouville measure on the unit tangent bundle $S\Sd$, which we identify with the set of orthonormal pairs $(x,v)$ in $\mathbb R^{d+1}$. The group $O(d+1)$ acts transitively on this set, and $\dd v\dd\mu(x)$ is an $O(d+1)$-invariant probability measure; hence it is the unique one. The geodesic flow is
\[
  \varphi_s(x,v)=\bigl(\cos s\,x+\sin s\,v,\ -\sin s\,x+\cos s\,v\bigr),
\]
that is, right multiplication of the $(d+1)\times2$ matrix $[\,x\ \ v\,]$ by a rotation $R_s\in SO(2)$. Right multiplication commutes with the left action of $O(d+1)$, so $\varphi_s$ pushes the invariant probability measure forward to an invariant probability measure, which by uniqueness is the same measure. In the same way, the push-forward of $\dd v\dd\mu(x)$ under the base-point projection $(x,v)\mapsto x$ is an $O(d+1)$-invariant probability measure on $\Sd$, hence equals $\mu$. Consequently
\[
  \iint g(\gamma_{x,v}(s))\dd v\dd\mu(x)=\int g\dd\mu
\]
for every $s$ and every measurable $g\ge0$. Taking $g=|\nabla f|$, the inner double integral equals $\int|\nabla f|\dd\mu$ for every $s$, and the above term in the bracket is $\theta\int|\nabla f|\dd\mu$. Therefore
\[
  \|f-V_nf\|_1\le\Bigl(\omega_d\int_0^\pi|\Phi(\cos\theta)|\,\theta\sin^{d-1}\theta\dd\theta\Bigr)
  \int|\nabla f|\dd\mu .
\]
By Lemma~\ref{lem:loc} with $k=d+2$ and $\sin\theta\le\theta$, the above term in the bracket is at most
\[
  C\,n^d\int_0^\infty\theta^{d}(1+n\theta)^{-k}\dd\theta
  =\frac{C}{n}\int_0^\infty s^d(1+s)^{-k}\dd s=\frac{C_J}{n}. 
\]
\end{proof}

\subsection{Upper Regular Point Sets and Sampling Bounds}\label{subsec:upperreg}

Recall that a finite multiset $Y=\{y_1,\dots,y_M\}\subset\Sd$ is \emph{upper regular at scale $t_1$ with constant $K\ge1$} if every spherical cap of radius $r>0$ contains at most $K(1+rt_1)^d$ points of $Y$, counted with multiplicity. It is easy to show that well-separated point sets are automatically upper regular.

\begin{lemma}[Well-separation implies upper regular]\label{lem:sep}
Let $Y=\{y_1,\ldots, y_M\}\subset \Sd$ be a point set. If $\dist(y_i,y_j)\ge c/t_1$ for all $i\ne j,1\le i,j\le M$, then $Y$ is upper regular at scale
$t_1$ with $K=K(c,d)$.
\end{lemma}

\begin{proof}
There are $0<a_d\le b_d$ with $a_d\rho^d\le\mu(C_\rho)\le b_d\rho^d$ for any cap $C_\rho$ of radius $0<\rho\le\pi$. The caps of radius $\rho_0=\min(c/(2t_1),\pi)$ about the points in $Y$ are pairwise disjoint. If $r\le\pi$, those centered in a given cap of radius $r$ lie in a cap of radius $\min(r+\rho_0,\pi)$, so their number is at most $b_d(r+\rho_0)^d/(a_d\rho_0^d)\le K(1+rt_1)^d$. If $r>\pi$, the number of points is at most $1/(a_d\rho_0^d)\le K(1+rt_1)^d$ as well. 
\end{proof}

The following Lemma shows that the sampling bound of a polynomial on an upper regular point set is controlled by its $L_1$-norm. 

\begin{lemma}[Sampling bound]\label{lem:sampling}
Let $Y=\{y_1,\ldots, y_M\}\subset \Sd$ be upper regular at scale $t_1$ with a constant $K\ge1$. Then for every $n\ge t_1$ and every polynomial $R$ of degree $\le n$, it holds that
\[
  \sum_{j=1}^M|R(y_j)|\le A\,n^d\,\|R\|_1
\]
for some constant  $A=A(d,K)$ depending on $d$ and $K$.
\end{lemma}

\begin{proof}
Since $V^{h_2}_nR=R$ for polynomial $R$ of degree at most $n$, Lemma~\ref{lem:loc} with $k=d+1$ gives
\[
  |R(y)|\le C\,n^d\int_{\Sd}(1+n\dist(y,z))^{-k}|R(z)|\dd\mu(z).
\]
Summing over $j$ and using Tonelli, it suffices to show $S(z):=\sum_j(1+n\dist(y_j,z))^{-k}\le B$ for every $z$, with $B=B(d,K)$. Split the points into the group with $\dist(y_j,z)<1/n$ and the groups $G_i$,
$i\ge1$, with $2^{i-1}/n\le\dist(y_j,z)<2^i/n$. By the upper regularity property of $Y$ and $t_1\le n$, the first group has at most $K(1+t_1/n)^d\le K2^d$ points, each of weight $\le1$. The group $G_i$ lies in a cap of radius $2^i/n$, so it has at most $K(1+2^it_1/n)^d\le K(1+2^i)^d$ points, each of weight $\le(1+2^{i-1})^{-k}$. Hence
\[
  S(z)\le K2^d+K\sum_{i\ge1}(1+2^i)^d(1+2^{i-1})^{-k}=:B<\infty
\]
because $k>d$. Take $A=CB$. 
\end{proof}

\begin{remark}\label{rem:Ht}
In a single application of Lemma~\ref{lem:sampling} at degree $n=t$, only the count of points in caps of radius $\ge1/t$ is used. Covering a cap of radius $r$ by $O((rt)^d)$ caps of radius $1/t$ shows that the assumption of upper regularity at scale $t_1$ may be replaced there by: every cap of radius $1/t$ contains at most $K'$ points. The dyadic argument below applies the lemma at all degrees $n\in[t_1,2t]$ and therefore uses the upper regularity as stated. Since the number of points in a cap is monotone in its radius, the upper regularity at scale $t_1$ is equivalent, up to the value of $K$, to the single condition that every cap of radius $1/t_1$ contains at most $K'$ points.
\end{remark}

\section{Nested Spherical Designs from Upper Regular Point Sets with Optimal Order}\label{sec:upperregular}
In this section, we prove Theorem~\ref{thm:extension} by estimating $\Lambda(Y,t)$ for upper regular point sets using a dyadic estimation technique and show that such an estimate is sharp. A chain of nested spherical designs can then be derived based on Theorem~\ref{thm:extension}.

\subsection{Extension from Upper Regular Point Sets}\label{subsec:extension}

We now estimate $\Lam(Y,t)$ for upper regular spherical designs and prove Theorem~\ref{thm:extension} with the master criterion (Proposition~\ref{prop:master}). We use two observations to replace the Cauchy--Schwarz estimate in \cite{ZZ}. First, if $n\le t_1$ then $V_nP$ has degree $<t_1$, so the design property gives $\sum_jV_nP(y_j)=M\int V_nP\dd\mu=0$: the spherical design annihilates the low frequencies of $P$ \emph{exactly}. Second, the high-frequency remainder $P-V_nP$ is small in $L^1$-norm,
$\|P-V_nP\|_1\lesssim n^{-1}\int|\nabla P|\dd\mu$ (by the $L^1$ Jackson-type estimate), and a
well-spread $Y$ cannot over-count a polynomial of degree $\le n$: $\sum_j|R(y_j)|\lesssim n^d\|R\|_1$ (a sampling inequality derived from kernel localization). Splitting $P-V_{t_1}P$ into dyadic frequency blocks, the block at frequency $n$ contributes $\lesssim n^{d-1}$, and the sum over blocks is dominated by $n\sim t$ when $d\ge2$.

\begin{proposition}[Dyadic estimate]\label{prop:dyadic}
Let $d\ge1$, $1\le t_1<t$, and $Y=\{y_1,\ldots,y_M\}\subset \Sd$ be a spherical $t_1$-design that is upper regular at scale $t_1$ with a constant $K\ge1$. Then for every $P\in\partial\Omega$,  it holds that
\[
  \Bigl|\sum_{j=1}^MP(y_j)\Bigr|\le
  \begin{cases}
    A_*\,t^{d-1}, & d\ge2,\\[2pt]
    A_*\bigl(1+\log_2(2t/t_1)\bigr), & d=1,
  \end{cases}
  \qquad A_*=A_*(d,K).
\]
\end{proposition}

\begin{proof}
Put $n_k=2^kt_1$, $\tilde V_k:=V_{n_k}^{h_1}$, and let $k_*=\min\{k:\ n_k\ge2t\}$. Since $t_1<t$ we have $n_0<2t$, so $k_*\ge1$ and $n_{k_*-1}<2t$. Because $\deg P\le t\le n_{k_*}/2$ and $h_1\equiv1$ on $[0,\frac12]$, we have
$\tilde V_{k_*}P=P$. Telescoping $P$, we have
\[
  P=\tilde V_0P+\sum_{k=0}^{k_*-1}R_k,\qquad R_k:=\tilde V_{k+1}P-\tilde V_kP .
\]

\emph{Low frequencies.} $\tilde V_0P$ has degree $<n_0=t_1$ and $\int \tilde V_0P\dd\mu=h_1(0)\int P\dd\mu=0$. Since $Y$ is a spherical $t_1$-design, $\sum_j\tilde V_0P(y_j)=M\int \tilde V_0P\dd\mu=0$. Hence $\sum_jP(y_j)=\sum_{k<k_*}\sum_jR_k(y_j)$.

\emph{Each dyadic block.} $R_k$ is a polynomial of degree $<n_{k+1}$, and
$n_{k+1}\ge t_1$. By Lemma~\ref{lem:jackson} and $\int|\nabla P|\dd\mu=1$, we have
\[
  \|R_k\|_1\le\|P-\tilde V_{k+1}P\|_1+\|P-\tilde V_kP\|_1\le C_J\Bigl(\frac1{n_{k+1}}+\frac1{n_k}\Bigr)\le\frac{2C_J}{n_k}.
\]
By Lemma~\ref{lem:sampling} at degree $n_{k+1}$,
\[
  \sum_j|R_k(y_j)|\le A\,n_{k+1}^d\|R_k\|_1\le 2^{d+1}AC_J\,n_k^{\,d-1}.
\]

\emph{Summation.} For $d\ge2$, $n_k^{d-1}=2^{(k-k_*+1)(d-1)}n_{k_*-1}^{d-1}$, so
\[
  \sum_{k<k_*}n_k^{d-1}\le n_{k_*-1}^{d-1}\sum_{i\ge0}2^{-i(d-1)}\le2(2t)^{d-1}.
\]
For $d=1$, the sum equals $k_*\le1+\log_2(2t/t_1)$. Combining the three displays
gives the claim with $A_*=2^{2d+1}AC_J$. 
\end{proof}

\begin{remark}[Relation to the one-step estimate]
Use the design property only for $\tilde V_0P$, so that $\sum_jP(y_j)=\sum_j(P-\tilde V_0P)(y_j)$, and apply
Lemma~\ref{lem:sampling} once, at degree $t$, to $R=P-\tilde V_0P$. Since $\|R\|_1\le C_J/t_1$ by
Lemma~\ref{lem:jackson}, this gives $|\sum_jP(y_j)|\le AC_J\,t^d/t_1$, the
\emph{one-step estimate}. It has the same order as Proposition~\ref{prop:dyadic} when $t_1$ is
comparable to $t$, and it is weaker when $t_1=o(t)$. The one-step argument combines the sampling
constant $t^d$ at the top degree with the $L^1$ bound $1/t_1$ at the lowest scale. The dyadic
splitting combines, for each block, a sampling constant $\lesssim n^d$ with an $L^1$ bound
$\lesssim1/n$ at the same frequency $n$, so each block contributes $n^{d-1}$.
\end{remark}

We are now ready to prove Theorem~\ref{thm:extension} that gives nested spherical designs from upper regular point sets.
\begin{proof}[Proof of Theorem~\ref{thm:extension}]
By Proposition~\ref{prop:dyadic}, $\Lam(Y,t)\le A_*t^{d-1}$ for $d\ge 2$. Apply Proposition~\ref{prop:master}: $Nc_d/t\ge A_*t^{d-1}$
is exactly $N\ge(A_*/c_d)t^d$. Hence, for every
\[
  N\ \ge\ \max\Bigl(C_{1,d},\ \frac{A_*}{c_d}\Bigr)\,t^d,
\]
there is $X\subset\Sd$ with $|X|=N$ such that $X\uplus Y$ is a spherical $t$-design, where $X$ can be taken of the form described in Proposition~\ref{prop:master}.
For the last assertion, $N=\lceil\max(C_{1,d},A_*/c_d)\rceil t^d$ and $M\le C_Mt_1^d\le C_Mt^d$.
\end{proof}

\begin{remark}[Dependence of the constants]
All constants depend on $d$, and we have not tracked this dependence. It enters through the localization constants $C_k(d,h)$, the volume and covering constants of caps, and the constants $C_{1,d},c_d$ of \cite{ZZ}. The quantity $C_M$ is not a constant of this kind but the ratio $M/t_1^d$ (or any larger number): the bounds depend on $M$ and $t_1$ only through $C_M$, and they are uniform in $t>t_1$.
\end{remark}

\begin{remark} The assumption on $M$ is
in fact automatic under the upper regularity assumption: a cap of radius $\pi$ is all of $\Sd$, so
$M\le K(1+\pi t_1)^d\le K(1+\pi)^dt_1^d$. In particular, no condition of the form $t_1\ge\varepsilon t$ is needed. 
\end{remark}

\begin{remark}[The case $d=1$]
For $d=1$ the same holds with $N\ge\max\bigl(C_{1,d},A_*(1+\log_2(2t/t_1))/c_d\bigr)\,t$, where $d=1$.
\end{remark}

\begin{remark} What is and is not claimed.
\begin{enumerate}
  \item[(a)] Theorem~\ref{thm:ZZ1} and Lemma~\ref{lem:loc} are quoted. The constant $c_d$ must be uniform in $N\ge C_{1,d}t^d$; this is how it is stated in \cite[Thm.~1]{ZZ}, where, in the notation of \cite{ZZ}, $c_d=C_{2,d}=r_d/(36\sqrt d)$ depends only on $d$.
  
  \item[(b)] Upper regularity is essential for Lemma~\ref{lem:sampling}, but it is automatic for spherical $t_1$-designs with $M\lesssim t_1^d$ points (Proposition~\ref{prop:auto}); Theorem~\ref{thm:arbitrary} removes it from the statement.
  
  \item[(c)] The identical-constants part of the Zheng--Zhuang conjecture is not addressed.
\end{enumerate}
\end{remark}

\subsection{Sharpness of \texorpdfstring{$\Lambda(Y,t)$}{Lambda(Y,t)}}\label{subsec:sharpness}

Proposition~\ref{prop:dyadic} shows that $\Lam(Y,t)\le A_*t^{d-1}$ for upper regular point sets for $d\ge2$. We now show that $t^{d-1}$ is the correct order of $\Lam(Y,t)$ for every point set that is not
too large, so that the estimate cannot be improved, regardless of the structure of $Y$.

\begin{lemma}[Bump polynomials]\label{lem:bump}
There are constants $n_0\ge2$ and $C_b\ge1$, depending only on $d$ and $h_1$, such
that for every $n\ge n_0$ and every $y\in\Sd$ there is a polynomial $Q_b=Q_{n,y}$ on
$\Sd$ of degree $\le2n-2$ with
\begin{enumerate}
  \item[(i)] $Q_b\ge0$ and $Q_b(y)=1$;
  \item[(ii)] $\displaystyle\int_{\Sd}Q_b\dd\mu\le C_b\,n^{-d}$;
  \item[(iii)] $\displaystyle\int_{\Sd}|\nabla Q_b|\dd\mu\le C_b\,n^{1-d}$.
\end{enumerate}
\end{lemma}

\begin{proof}
Write $Z_\ell(\langle x,y\rangle)=D(\ell,d)P_{\ell,d}(\langle x,y\rangle)$, the
reproducing kernel of the degree-$\ell$ spherical harmonics, and set
\[
  g(x)=\sum_{\ell<n}h_1(\ell/n)\,Z_\ell(\langle x,y\rangle)=\Phi^{h_1}_n(\langle x,y\rangle).
\]
Since $h_1(\ell/n)=0$ for $\ell\ge n$, $g$ is a polynomial of degree $\le n-1$. Define
\[
  Q_b=\frac{g^2}{g(y)^2}.
\]
Then $Q_b\ge0$ is a polynomial of degree $\le2n-2$, and $Q_b(y)=1$, provided $g(y)\ne0$.

\emph{Size of $g(y)$.} Because $P_{\ell,d}(1)=1$,
$g(y)=\sum_{\ell<n}h_1(\ell/n)D(\ell,d)$. All terms are $\ge0$ and $h_1=1$ on
$[0,\frac12]$, so
\[
  g(y)\ \ge\ \sum_{\ell\le n/2}D(\ell,d)=\dim\Pi_{\lfloor n/2\rfloor}(\Sd)\ \ge\ a'_d\,n^d
\]
for $n\ge2$, with $a'_d>0$ depending only on $d$ (the dimension of the space of
polynomials of degree $\le m$ on $\Sd$, constants included, is comparable to $m^d$, and $\lfloor n/2\rfloor\ge n/3$ for $n\ge2$). In particular
$g(y)\neq0$.

\emph{Spectral identities.} The functions $Z_\ell(\langle\cdot,y\rangle)$ lie in
mutually orthogonal eigenspaces of the Laplace--Beltrami operator, with eigenvalue
$-\ell(\ell+d-1)$, and $\int Z_\ell(\langle x,y\rangle)^2\dd\mu(x)=D(\ell,d)$. Hence
\[
  \|g\|_2^2=\sum_{\ell<n}h_1(\ell/n)^2D(\ell,d)\le\sum_{\ell<n}D(\ell,d)\le b'_d\,n^d,
\]
\[
  \|\nabla g\|_2^2=\int g\,(-\Delta g)\dd\mu=\sum_{\ell<n}\ell(\ell+d-1)\,h_1(\ell/n)^2D(\ell,d)\le(d+1)\,n^2\,b'_d\,n^d ,
\]
where $b'_d$ bounds $\dim\Pi_{n-1}(\Sd)/n^d$ and we used
$\ell(\ell+d-1)\le n(n+d)\le(d+1)n^2$ for $\ell<n$.

\emph{(ii).} $\int Q_b\dd\mu=\|g\|_2^2/g(y)^2\le b'_dn^d/(a_d'^{\,2}n^{2d})=C\,n^{-d}$.

\emph{(iii).} Since $|\nabla Q_b|=2|g||\nabla g|/g(y)^2$, Cauchy--Schwarz gives
\[
  \int|\nabla Q_b|\dd\mu\le\frac{2\|g\|_2\|\nabla g\|_2}{g(y)^2}
  \le\frac{2\,(b'_dn^d)^{1/2}\,\bigl((d+1)b'_d\bigr)^{1/2}n^{(d+2)/2}}{a_d'^{\,2}\,n^{2d}}=C\,n^{1-d}.
\]
Take $C_b$ as the largest of the constants in (ii) and (iii) and $n_0$ as given. \qedhere
\end{proof}

\begin{remark}[Numerical illustration, not part of the proof]
Take $d=2$ and, for $s\in[\tfrac12,1]$, $h_1(s)=\psi(1-s)/\bigl(\psi(1-s)+\psi(s-\tfrac12)\bigr)$,
where $\psi(x)=e^{-1/x}$ for $x>0$ and $\psi(x)=0$ otherwise, with $h_1=1$ on $[0,\tfrac12]$ and
$h_1=0$ on $[1,\infty)$. Computation with Legendre expansions and numerical quadrature gives the
following values.
\begin{center}
\begin{tabular}{lccccc}
\hline
$n$ & 8 & 16 & 32 & 64 & 128\\
\hline
$n^2\int Q_b\dd\mu$ & 1.39 & 1.51 & 1.57 & 1.61 & 1.62\\
$n\int|\nabla Q_b|\dd\mu$ & 1.14 & 1.21 & 1.24 & 1.25 & 1.26\\
\hline
\end{tabular}
\end{center}
Both quantities appear to level off, consistent with Lemma~\ref{lem:bump}(ii) and (iii). The limiting values depend on the choice of $h_1$.
\end{remark}

\begin{theorem}[Sharpness]\label{thm:sharpness}
There are constants $c_0,c_1>0$ and $t_0\in \mathbb{N}_0$, depending only on $d$, such that if $t\ge t_0$ and $1\le M\le c_0\,t^d$ then $\Lam(Y,t)\ge c_1\,t^{d-1}$ for every multiset $Y$ of $M$ points. 
\end{theorem}
\begin{proof}
Let $n=\lfloor t/2\rfloor+1$, so that $2n-2\le t$ and $n\ge t/2$, and assume $t\ge t_0$
with $t_0$ large enough that $n\ge n_0$ and $C_bn^{-d}\le\frac12$.
Pick $y_1\in Y$ and let $Q_b=Q_{n,y_1}$ be the bump polynomial of Lemma~\ref{lem:bump} at $y_1$. Put
\[
  \beta=\int|\nabla Q_b|\dd\mu,\qquad \bar Q_b=\int Q_b\dd\mu,\qquad
  P_0=\frac{Q_b-\bar Q_b}{\beta}.
\]
$Q_b$ is not constant, because $Q_b(y_1)=1$ while $\bar Q_b\le C_bn^{-d}\le\frac12$; hence
$\beta>0$. $P_0$ is a polynomial of degree $\le t$ with mean zero and
$\int|\nabla P_0|\dd\mu=1$, that is, $P_0\in\partial\Omega$. Using $Q_b\ge0$ and
$y_1\in Y$, we have
\[
  \sum_{y\in Y}P_0(y)=\frac{1}{\beta}\Bigl(\sum_{y\in Y}Q_b(y)-M\bar Q_b\Bigr)
  \ge\frac1\beta\bigl(Q_b(y_1)-M\bar Q_b\bigr)\ge\frac1\beta\bigl(1-MC_bn^{-d}\bigr).
\]
If $M\le n^d/(2C_b)$, which holds when $M\le c_0t^d$ with $c_0=2^{-d-1}/C_b$ since
$n\ge t/2$, the bracket is $\ge\frac12$, and by Lemma~\ref{lem:bump}(iii),
$\beta\le C_bn^{1-d}\le C_b\,2^{d-1}t^{1-d}$ for $d\ge1$ (using $n\ge t/2$ when
$d\ge2$, and $n^{1-d}=1$ when $d=1$). Therefore,
\[
  \sum_{y\in Y}P_0(y)\ \ge\ \frac{t^{d-1}}{2^{d}\,C_b}.
\]
Since $\partial\Omega=-\partial\Omega$, the polynomial $-P_0\in\partial\Omega$ gives
$\Lam(Y,t)\ge\sum_yP_0(y)\ge c_1t^{d-1}$ with $c_1=2^{-d}C_b^{-1}$. \qedhere
\end{proof}

Using the above result, we can show that the order of $\Lambda(Y,t)$ is exactly $t^{d-1}$ for an upper regular spherical $t_1$-design $Y$ with at most $c_0t^d$ points.

\begin{corollary}[Order of $\Lam$ for upper regular spherical designs]\label{cor:order}
Let $c_0,c_1,t_0$ be the constants of Theorem~\ref{thm:sharpness}, let $d\ge2$ and $t\ge t_0$.
Let $Y=\{y_1,\ldots,y_M\}\subset\Sd$ be a spherical $t_1$-design with $1\le t_1<t$ and $M\le c_0t^d$ points
that is upper regular at scale $t_1$ with a constant $K\ge1$. This holds, for instance, if the points of $Y$
have minimal separation distance at least $c/t_1$, with $K=K(c,d)$ by Lemma~\ref{lem:sep}. Then
\[
  c_1\,t^{d-1}\ \le\ \Lam(Y,t)\ \le\ A_*(d,K)\,t^{d-1}.
\]
\end{corollary}

\begin{proof}
The upper bound is Proposition~\ref{prop:dyadic}. The lower bound is Theorem~\ref{thm:sharpness}, which applies because $t\ge t_0$ and $1\le M\le c_0t^d$. The size condition is not restrictive when $t_1\ll t$: a cap of radius $\pi$ is all of $\Sd$, so $M\le K(1+\pi t_1)^d\le K(1+\pi)^dt_1^d$, and hence $M\le c_0t^d$ holds as soon as $t_1\le\bigl(c_0/(K(1+\pi)^d)\bigr)^{1/d}\,t$. \qedhere
\end{proof}

\begin{remark}
Consequently, when $t_1=o(t)$ the factor $1/t_1$ in the one-step bound $\Lam\lesssim t^d/t_1$ is an artifact of that estimate, since the order of $\Lam$ does not depend on $t_1$.
\end{remark}

\begin{remark}\label{rem:meaning} We have the following remarks on the sharpness of $\Lambda$.
\begin{itemize}
\item[(a)] Theorem~\ref{thm:sharpness} pins down what the Brouwer-degree scheme can deliver.
By Proposition~\ref{prop:master}, the conditions $N\ge C_{1,d}t^d$ and
$Nc_d/t\ge\Lam(Y,t)$ suffice to complete $Y$, and Theorem~\ref{thm:sharpness} shows
that the second forces $N\ge c_1t^d/c_d$ whenever $1\le M\le c_0t^d$. With the
map of Theorem~1 of Zheng--Zhuang, then, the scheme cannot certify fewer than
order $t^d$ new points. This agrees with the classical bound $N+M\ge ct^d$ for
$t$-designs and adds nothing to it. Since $\Lam$ is only a sufficient quantity,
some other method might complete $Y$ with fewer points, and we prove no lower
bound for the nested problem beyond the classical one.

\item[(b)] The size restriction $M\le c_0t^d$ is used. For $M\gg t^d$ the correction
$M\bar Q_b$ can dominate $Q_b(y_1)$, and a bump at one point no longer forces a
large negative value. A lower bound of the same order in that range would need
a different test function.

\item[(c)] The theorem also sharpens the picture of $\Lam$ itself. The $L^2$ estimate of
Zheng--Zhuang gives $\Lam\lesssim Mt^d$ and the one-step estimate of
Section~\ref{subsec:extension} gives $\Lam\lesssim t^d/t_1$, while
Corollary~\ref{cor:order} gives $\Lam\asymp t^{d-1}$ for upper regular spherical designs.
The first bound is therefore far from the truth, and the second is far from it
once $t_1\ll t$. Coincident points show the same gap. If $y_1$ has multiplicity
$m$ in $Y$ and $M\le c_0t^d$, the proof of Theorem~\ref{thm:sharpness} yields
$\Lam\ge c_1\,m\,t^{d-1}$, since the bracket there is at least $m-\tfrac12\ge m/2$.
For $m=M$ this meets the upper bound $CMt^{d-1}$ of Remark~\ref{rem:SN}, valid
for $d\ge2$, and the $L^2$ bound overshoots it by a factor of $t$.

\item[(d)] When $d=1$ we have $t^{d-1}=1$, so Theorem~\ref{thm:sharpness} gives only $\Lam\ge c_1$,
while the dyadic estimate (Proposition~\ref{prop:dyadic}) gives
$\Lam\lesssim1+\log(t/t_1)$. For the regular polygon,
Section~\ref{subsec:circle} shows $\Lam\le\pi$ for all $t$. For general spherical designs
on the circle the question is open.
\end{itemize}
\end{remark}

\subsection{Chain of Nested Spherical Designs}\label{subsec:chains}

We now iterate Theorem~\ref{thm:extension} to give chains of nested spherical designs. The point to check is that upper regularity
survives the extension step. The relevant property of the new points is not in the
statement of Theorem~\ref{thm:ZZ1} but in its proof, which we recall.

\begin{lemma}[Localization of the new points; properties of the construction in {\cite[Appendix~B]{ZZ}}]\label{lem:local}
Let $t\ge1$, let $N=\lambda t^d$ be an integer with $\lambda\ge C_{1,d}$ (enlarging $C_{1,d}$ if necessary, we assume $C_{1,d}\ge1$), and let $F$ be the map of Theorem~\ref{thm:ZZ1} as constructed in \cite[Appendix~B]{ZZ}. There are constants $b_d,B_d,r_d>0$ depending only on $d$ such that the measure of any spherical cap $C(z,\rho)$ centered at $z$ with radius $\rho$ satisfies  $\mu(C(z,\rho))\le b_d\rho^d$ for all $\rho>0$, and for every $Q\in\Pt$ the point set $X=\{x_i(Q)\}_{i=1}^N$ has the following property: there is an area-regular partition $\{R_i\}_{i=1}^N$ of $\Sd$ (independent
of $Q$) with $\operatorname{diam}R_i\le B_dN^{-1/d}$, and points $\tilde x_i\in R_i$ with
$\dist(x_i(Q),\tilde x_i)\le r_d/(3t)$. Consequently, every cap of radius $r>0$ contains
at most $K_X(1+rt)^d$ points of $X$, where
\[
  K_X=\beta_d\,\lambda,\qquad\beta_d=b_d\max\bigl(1,\tfrac{r_d}{3}+B_d\bigr)^d .
\]
\end{lemma}

\begin{proof} In the construction in \cite[Appendix~B, (B.1)--(B.3)]{ZZ},
$x_i(P)=w_i(P,r_d/(3t))$, where $w_i$ solves $w'=U(P,w)$ with
$w(0)=\tilde x_i\in R_i$ and $|U|\le1$; hence $\dist(x_i(P),\tilde x_i)\le r_d/(3t)$. The
partition is area-regular with norm $\le B_dN^{-1/d}$. Suppose $x_i(Q)\in C(z,r)$. Then
$\tilde x_i\in C(z,r+r_d/(3t))$ and $R_i\subset C(z,\rho)$ with
$\rho=r+r_d/(3t)+B_dN^{-1/d}\le r+c'/t$, where $c'=r_d/3+B_d$ and we used
$N^{-1/d}=\lambda^{-1/d}/t\le1/t$. The cells have area $1/N$ and overlap in null sets, so
the number of such $i$ is at most $N\mu(C(z,\rho))\le b_d\lambda t^d(r+c'/t)^d
=b_d\lambda(rt+c')^d\le K_X(1+rt)^d$. \qedhere
\end{proof}

\begin{remark}
The new points do \emph{not} lie in the cells themselves; they are displaced by up to
$r_d/(3t)$ by the flow. This is harmless here, because the displacement is on the scale
$1/t$ at which upper regularity is measured.
\end{remark}

\begin{lemma}[Propagation of upper regularity] \label{lem:propagate}
Let $Y$ be upper regular at scale $s$ with a constant $K$, let $s\le t$, and let $X,K_X$ be as
in Lemma~\ref{lem:local}. Then $X\uplus Y$ is upper regular at scale $t$ with
a constant $K+K_X$.
\end{lemma}
\begin{proof}
A cap of radius $r$ contains at most $K(1+rs)^d\le K(1+rt)^d$ points of $Y$ and at most
$K_X(1+rt)^d$ points of $X$. \qedhere
\end{proof}

By inspection of Lemma~\ref{lem:sampling} and Proposition~\ref{prop:dyadic}, the constant
$A_*(d,K)$ of Proposition~\ref{prop:dyadic} is linear in $K$:
$A_*(d,K)=\alpha_d K$ with $\alpha_d$ depending only on $d$.

\begin{theorem}[Finite chains]\label{thm:chain}
Let $d\ge2$ and $L\ge2$, let $1\le t_1<t_2<\dots<t_L$ be integers, and let $Y_0$ be a spherical
$t_1$-design that is upper regular at scale $t_1$ with constant $K_0$ and $|Y_0|\le C_Mt_1^d$, $C_M\ge1$.
Then, there are multisets $X_1,\dots,X_{L-1}$ such that $Y_k=Y_{k-1}\uplus X_k$ is a spherical
$t_{k+1}$-design for $k=1,\dots,L-1$, with
\[
  |X_k|\le C(d,K_0)\,\gamma_d^{\,k}\,t_{k+1}^d,\qquad
  |Y_{L-1}|\le C(d,K_0,C_M)\,\gamma_d^{\,L}\,t_L^d,
\]
where $\gamma_d>1$ depends only on $d$. 
\end{theorem}
\begin{proof}
Let $K_k$ be a constant for which $Y_k$ is upper regular at scale $t_{k+1}$. Suppose
this holds for $Y_{k-1}$ at scale $t_k$. Since $Y_{k-1}$ is a spherical $t_k$-design and $t_k<t_{k+1}$,
Theorem~\ref{thm:extension} gives $X_k$ with $|X_k|=N_k=\lceil\lambda_k\rceil t_{k+1}^d$, where
$\lambda_k=\max(C_{1,d},\alpha_dK_{k-1}/c_d)\ge1$, such that $Y_k$ is a spherical $t_{k+1}$-design.
Moreover, by the last assertion of Proposition~\ref{prop:master}, $X_k$ has the form of Lemma~\ref{lem:local} (with $\lceil\lambda_k\rceil\le2\lambda_k$
in place of $\lambda$), so by Lemma~\ref{lem:propagate} $Y_k$ is upper regular at scale
$t_{k+1}$ with
\[
  K_k\le K_{k-1}+2\beta_d\lambda_k\le\gamma_dK_{k-1}+\delta_d,\qquad
  \gamma_d=1+\frac{2\beta_d\alpha_d}{c_d},\ \ \delta_d=2\beta_dC_{1,d}.
\]
Hence $K_k+\delta_d/(\gamma_d-1)\le\gamma_d^{\,k}(K_0+\delta_d/(\gamma_d-1))$, which bounds
$\lambda_k$ and therefore $|X_k|$. Summing, $|Y_{L-1}|\le C_Mt_1^d+\sum_{k=1}^{L-1}|X_k|
\le C\gamma_d^{\,L}t_L^d$, using $t_{k+1}\le t_L$. 
\end{proof}

\begin{remark}\label{rem:chain} We have the following remarks.
\begin{itemize}
\item[(i)] Note that no condition on the ratios $t_{k+1}/t_k$ is needed in the above result. The independence of the ratios $t_{k+1}/t_k$ is a consequence of Theorem~\ref{thm:extension}
being uniform in $t_1<t$. With only the one-step estimate of Section~\ref{subsec:extension} (of order $t^d/t_1$) one would
need a bound $t_{k+1}\le mt_k$ to keep each level of optimal order.

\item[(ii)] The growth $\gamma_d^{\,L}$ comes from two sources: $K_X$ is proportional to $\lambda$, and
$\lambda$ must be proportional to $K_{k-1}$ for the master criterion to hold. For a fixed
number of levels the result is of optimal order; for $L\to\infty$ (for instance
$t_{k+1}=2t_k$ with $L\sim\log t$) it gives only $t^{d+\log_2\gamma_d}$. Whether
the constants can be made independent of $L$ is open.

\item[(iii)] Zheng and Zhuang in \cite{ZZ} remarked that one may seek multi-nested structures such as
$t_1=2t_2=4t_3=\cdots$ with identical constants \cite[Section~5]{ZZ}; their
Proposition~1 produces these trivially by reusing a spherical $t$-design, whereas
Theorem~\ref{thm:chain} starts from an arbitrary upper regular spherical design (and, by Proposition~\ref{prop:auto}, from any spherical design with $M\lesssim t_1^d$ points).
\end{itemize}
\end{remark}

\section{Nested Spherical Designs from Arbitrary Spherical Designs}\label{sec:beyond}
In this section, we go beyond upper regular point sets. Section~\ref{subsec:arbitrary} shows that the dyadic argument depends on the starting point set only through its $L^1$ sampling constant (Theorem~\ref{thm:master}), relates this constant to the local multiplicity of the point set (Proposition~\ref{prop:mult}), and proves that bounded multiplicity at scale $1/t$ is necessary for any extension of optimal order (Proposition~\ref{prop:necessity}). Applying the necessity statement to a spherical $t_1$-design itself shows that every spherical $t_1$-design with $M\lesssim t_1^d$ points is automatically upper regular (Proposition~\ref{prop:auto}), which yields Theorem~\ref{thm:arbitrary}. Section~\ref{subsec:circle} treats the circle $d=1$, where a direct discrepancy argument removes the logarithmic factor (Theorem~\ref{thm:circle}), and compares the resulting bounds with numerically computed values.

\subsection{Extension from Arbitrary Point Sets}\label{subsec:arbitrary}
 
The dyadic argument of Section~\ref{subsec:extension} uses upper regularity of the starting point set only through one quantity, the $L^1$ sampling constant. We isolate it, record what it gives for arbitrary spherical designs, and show that bounded multiplicity at scale $1/t$ cannot be dispensed with.
 
\begin{definition}\label{def:sampling}
For a finite multiset $Y=\{y_1,\dots,y_M\}\subset\Sd$ and $n\ge1$ the \emph{$L^1$
sampling constant} of $Y$ is
\[
  A_Y(n)=\sup\Bigl\{\frac{\sum_{j=1}^M|R(y_j)|}{\|R\|_1}:\ R\ne0\text{ a polynomial of degree}\le n\Bigr\}.
\]
\end{definition}

\begin{theorem}[Master estimate]\label{thm:master}
Let $1\le t_1<t$ and let $Y\subset \Sd$ be a spherical $t_1$-design. With $n_k=2^kt_1$ and
$k_*=\min\{k:n_k\ge2t\}$,
\[
  \Lam(Y,t)\ \le\ 2C_J\sum_{k=0}^{k_*-1}\frac{A_Y(n_{k+1})}{n_k},
\]
where $C_J$ is the constant of Lemma~\ref{lem:jackson}.
\end{theorem}
\begin{proof}
Repeat the proof of Proposition~\ref{prop:dyadic}: for $P\in\partial\Omega$ write
$P=\tilde V_0P+\sum_{k<k_*}R_k$ with $\tilde V_k=V^{h_1}_{n_k}$; the low-frequency part contributes
nothing because $Y$ is a spherical $t_1$-design, $\|R_k\|_1\le2C_J/n_k$, and $R_k$ has degree $<n_{k+1}$.
Instead of Lemma~\ref{lem:sampling}, by the definition of $A_Y(n)$ in Definition~\ref{def:sampling}, we have $\sum_j|R_k(y_j)|\le A_Y(n_{k+1})\|R_k\|_1$. Then, summing over $k$ gives the conclusion. 
\end{proof}

Two cases are already covered by the master estimate. If $A_Y(n)\lesssim n^d$ for $n\ge t_1$, which holds for upper regular point sets by Lemma~\ref{lem:sampling}, then Theorem~\ref{thm:master} recovers Proposition~\ref{prop:dyadic} and hence Theorem~\ref{thm:extension}. For an arbitrary spherical $t_1$-design $Y$ with $M$ points one has $A_Y(n)\lesssim Mn^d$, and Theorem~\ref{thm:master} gives $\Lam(Y,t)\lesssim Mt^{d-1}$, so that $N\ge CMt^d$ new points suffice, one power of $t$ better than \cite{ZZ}. The second case is summarized in the following corollary.


\begin{corollary}[Arbitrary spherical designs: one power of $t$ better than \cite{ZZ}]\label{cor:arbitrary}
Let $d\ge2$, $1\le t_1<t$, and let $Y$ be an arbitrary spherical $t_1$-design with $M$
points. Then $\Lam(Y,t)\le C\,M\,t^{d-1}$ with $C=C(d)$, so every
$N\ge\max(C_{1,d},\,C M/c_d)\,t^d$ admits $X$ with $|X|=N$ such that $X\uplus Y$ is a
spherical $t$-design. In particular, if $M\le C_Mt_1^d$, then $N+M\le C(d,C_M)\,t^{2d}$.
\end{corollary}
\begin{proof}
For a polynomial $R$ of degree $\le n$, $V^{h_2}_nR=R$, so Lemma~\ref{lem:loc} with $k=1$
gives $|R(y)|\le C\,n^d\|R\|_1$ for every $y$. Hence $A_Y(n)\le C\,M\,n^d$, and
$\sum_{k<k_*}n_{k+1}^d/n_k=2^d\sum_{k<k_*}n_k^{d-1}\le2^{d+1}(2t)^{d-1}$ for $d\ge2$.
Apply Theorem~\ref{thm:master} and Proposition~\ref{prop:master}. \qedhere
\end{proof}
 
\begin{remark}\label{rem:SN}
For comparison, \cite[Thm.~3 and Cor.~2]{ZZ} give $N\gtrsim Mt^{d+1}$ for arbitrary spherical $t_1$-designs
and, for $M\sim t_1^d$, a total of order $t^{2d+1}$. Corollary~\ref{cor:arbitrary}
improves this to $t^{2d}$ for every spherical $t_1$-design of optimal order, and Theorem~\ref{thm:arbitrary}
below improves it further to $t^d$. The spherical design property is not needed for the bound $\Lam\le CMt^{d-1}$. For an arbitrary multiset $Y$ of $M$ points, it follows
from the inequality $\|P\|_\infty\le Ct^{d-1}\int|\nabla P|\dd\mu$ on $\Pt$, which combines the
Sobolev inequality $\|P\|_{d/(d-1)}\lesssim\int|\nabla P|\dd\mu$ for mean-zero $P$ with the
Nikolskii inequality $\|P\|_\infty\lesssim t^{d-1}\|P\|_{d/(d-1)}$ ($d\ge2$, see \cite{DX13}).
Corollary~\ref{cor:arbitrary} is included because it is the special case
$A_Y(n)\lesssim Mn^d$ of Theorem~\ref{thm:master}.
\end{remark}

We next show that bounded multiplicity of $Y$ at scale $1/t$ is \emph{necessary} for a nested extension of optimal order (see Proposition~\ref{prop:necessity}), so some regularity of $Y$ cannot be dispensed with.

For a finite multiset $Y$ and $r>0$ let $\kappa_Y(r)=\max_{z\in\Sd}\#\bigl(Y\cap C(z,r)\bigr)$,
where $C(z,r)$ is the closed cap of radius $r$ about $z$ and points are counted with
multiplicity.
 
\begin{lemma}[Bernstein inequality on the sphere]\label{lem:bernstein}
If $P$ is a spherical polynomial of degree at most $n$, then $\|\nabla P\|_\infty\le n\|P\|_\infty$.
\end{lemma}
\begin{proof}
Fix $x\in\Sd$ and a unit tangent vector $v$ at $x$, and let $\gamma(s)=\cos s\,x+\sin s\,v$.
Since $P$ is the restriction of a polynomial of degree at most $n$ on $\mathbb R^{d+1}$ and the
coordinates of $\gamma$ are linear in $\cos s$ and $\sin s$, the function $T(s)=P(\gamma(s))$ is a
trigonometric polynomial of degree at most $n$. By the classical Bernstein inequality,
$|T'(0)|\le n\|T\|_\infty\le n\|P\|_\infty$. Since $T'(0)=\langle\nabla P(x),v\rangle$, choosing
$v=\nabla P(x)/|\nabla P(x)|$ gives $|\nabla P(x)|\le n\|P\|_\infty$. \qedhere
\end{proof}
 
\begin{proposition}\label{prop:mult}
There are constants $c_*\in(0,1]$, $C_*$ and $n_*\ge2$, depending only on $d$, such that
for every finite multiset $Y$ and every $n\ge n_*$,
\[
  c_*\,n^d\,\kappa_Y(c_*/n)\ \le\ A_Y(n)\ \le\ C_*\,n^d\sup_{z\in\Sd}\sum_{j}\bigl(1+n\dist(y_j,z)\bigr)^{-(d+1)} .
\]
\end{proposition}
\begin{proof}
The upper bound is the first display in the proof of Lemma~\ref{lem:sampling}. For the
lower bound let $m=\lfloor n/2\rfloor+1\le n$, so $2m-2\le n$, choose $z$ with
$\#(Y\cap C(z,c_*/n))=\kappa_Y(c_*/n)$, and let $Q_b=Q_{m,z}=g^2/g(z)^2$ be the bump polynomial of
Lemma~\ref{lem:bump}, which has degree $\le n$ and $\|Q_b\|_1=\int Q_b\dd\mu\le C_bm^{-d}$.
All coefficients of $g=\Phi^{h_1}_m(\langle\cdot,z\rangle)$ are nonnegative and
$|P_{\ell,d}|\le1$, so $|g|\le g(z)$. Since $g$ has degree at most $m-1$ and $|g|\le g(z)$, Lemma~\ref{lem:bernstein} gives
$\|\nabla g\|_\infty\le m\,g(z)$, hence $g(x)\ge g(z)\bigl(1-m\operatorname{dist}(x,z)\bigr)$ and
$Q_b(x)\ge\tfrac12$ whenever $\operatorname{dist}(x,z)\le c/m$, with $c=1-2^{-1/2}$. Since $m\le n$, the cap $C(z,c_*/n)$ lies in $C(z,c/m)$ as soon as
$c_*\le c$, so $\sum_j|Q_b(y_j)|\ge\frac12\kappa_Y(c_*/n)$. Moreover $\|Q_b\|_1\le C_bm^{-d}\le2^dC_bn^{-d}$
because $m>n/2$. Hence
$A_Y(n)\ge\sum_j|Q_b(y_j)|/\|Q_b\|_1\ge2^{-d-1}C_b^{-1}\,n^d\,\kappa_Y(c_*/n)$, and the claim holds with
$c_*=\min\bigl(1,c,2^{-d-1}C_b^{-1}\bigr)$ and $n_*=2n_0$ (so that $m\ge n_0$ in Lemma~\ref{lem:bump}). \qedhere
\end{proof}
 
Thus a cluster of $m$ points of $Y$ within distance $c_*/n$ forces $A_Y(n)\gtrsim m\,n^d$, while
for a point set that is upper regular at scale $t_1$ the right-hand side is $\lesssim n^d$ for $n\ge t_1$ and, for nonempty $Y$,
$\kappa_Y\ge1$, so both sides are of order $n^d$. In this sense upper regularity is the natural
way to keep the sampling constant of order $n^d$ for all $n\in[t_1,2t]$. Conversely, if
$t_1\ge n_*$ and $A_Y(t_1)\le A\,t_1^d$, the lower bound gives $\kappa_Y(c_*/t_1)\le A/c_*$, and
covering a cap of radius $r$ by $O((1+rt_1)^d)$ caps of radius $c_*/t_1$ shows that $Y$ is upper regular at
scale $t_1$ with a constant depending only on $A$ and $d$.
 
\begin{proposition}[Necessity]\label{prop:necessity}
There are constants $c_*,C_{**}$ and $t_{**}$, depending only on $d$, such that for $t\ge t_{**}$
the following holds. If $Z$ is a spherical $t$-design that contains $Y$ as a sub-multiset,
and $T=|Z|$, then
\[
  \kappa_Y(c_*/t)\ \le\ C_{**}\,T\,t^{-d}.
\]
In particular, if $Y$ extends to a spherical $t$-design with at most $\Gamma t^d$ points,
then every cap of radius $c_*/t$ contains at most $C_{**}\Gamma$ points of $Y$.
\end{proposition}
\begin{proof}
Let $m=\lfloor t/2\rfloor+1$, let $z$ be a point where $\kappa_Y(c_*/t)$ is attained, and let
$Q_b=Q_{m,z}$ be the bump polynomial of Lemma~\ref{lem:bump}, a nonnegative polynomial of degree $\le t$. Since $Z$ is a
spherical $t$-design, $\sum_{w\in Z}Q_b(w)=T\int Q_b\dd\mu\le TC_bm^{-d}$. On the other hand, as in the proof of Proposition~\ref{prop:mult},
$Q_b\ge\frac12$ on $C(z,c_*/t)$ and $Q_b\ge0$ elsewhere, so
$\sum_{w\in Z}Q_b(w)\ge\frac12\#(Y\cap C(z,c_*/t))$. Combine with $m\ge t/2$.
\end{proof}
 
\begin{remark} We have the following remarks.
\begin{itemize}
\item[(i)] If $Y$ consists of $m$ copies of a spherical $t_1$-design, then $Y$ is again a spherical $t_1$-design, and
Proposition~\ref{prop:necessity} shows that every spherical $t$-design containing it has at least
$c\,m\,t^d$ points. So some regularity of $Y$ is unavoidable, not merely a convenience
of the method; this is a lower bound for the nested problem itself, not only for the
Brouwer scheme.

\item[(ii)] Upper regularity at scale $t_1$ asks for more than Proposition~\ref{prop:necessity} makes
necessary. It is equivalent, up to constants, to bounded multiplicity at scale $1/t_1$
(Remark~\ref{rem:Ht}), whereas necessity is shown only at the smaller scale $1/t$. For spherical
$t_1$-designs with $M\lesssim t_1^d$ points the gap closes (Proposition~\ref{prop:auto} below). For a general point set
$Y$, whether bounded multiplicity at scale $1/t$ already suffices is open (Section~\ref{sec:conclusion}).  
\end{itemize}
\end{remark}

\begin{remark}[Upper regularity here and in Section~\ref{sec:upperregular}] \label{rem:regvs}
In Theorem~\ref{thm:extension}, upper regularity is a hypothesis on the starting set. It is a geometric condition on the number of points in caps, satisfied for instance by well-separated sets, while the design property is used only to cancel the low frequencies. In Proposition~\ref{prop:auto} the same notion is a consequence: for a spherical $t_1$-design it follows from the design property together with the size bound $M\le C_Mt_1^d$, by testing the design against a bump polynomial. Hence Theorem~\ref{thm:arbitrary} has no regularity hypothesis, although its proof passes through Theorem~\ref{thm:extension}. What goes beyond Section~\ref{sec:upperregular} is therefore the class of starting sets (every spherical $t_1$-design of optimal order, separated or not). For a point set that is not a spherical design no such statement holds: Proposition~\ref{prop:necessity} shows that bounded multiplicity at scale $1/t$ is necessary, and one has only the bound $\Lam\le CMt^{d-1}$ of Remark~\ref{rem:SN}.
\end{remark}

\begin{proposition}[Automatic regularity]\label{prop:auto}
Let $d\ge1$ and $C_M\ge1$. There is $K=K(d,C_M)$, linear in $C_M$, such that every spherical
$t_1$-design $Y$ with $M=|Y|\le C_Mt_1^d$ is upper regular at scale $t_1$ with the constant $K$.
\end{proposition}
\begin{proof}
Let $c_*,C_{**},t_{**}$ be as in Proposition~\ref{prop:necessity}. If $t_1\ge t_{**}$, apply that
proposition with $Z=Y$, $t=t_1$ and $T=M$. This gives
$\kappa_Y(c_*/t_1)\le C_{**}Mt_1^{-d}\le C_{**}C_M$. If $t_1<t_{**}$, then
$\kappa_Y(\cdot)\le M\le C_Mt_{**}^d$. In both cases
$\kappa_Y(c_*/t_1)\le\kappa_0:=\max(C_{**},t_{**}^d)\,C_M$. A cap of radius $r$ is covered by at most
$C_d(1+rt_1/c_*)^d$ caps of radius $c_*/t_1$, so it contains at most
$C_dc_*^{-d}\kappa_0(1+rt_1)^d$ points of $Y$. \qedhere
\end{proof}

Now we are ready to prove Theorem~\ref{thm:arbitrary}. 
\begin{proof}[Proof of Theorem~\ref{thm:arbitrary}]
By Proposition~\ref{prop:auto}, $Y$ is upper regular at scale $t_1$ with a constant $K=K(d,C_M)$, linear in $C_M$.
Apply Proposition~\ref{prop:dyadic} and Theorem~\ref{thm:extension}; since $A_*(d,K)=\alpha_dK$ and
$C(d,K)$ is linear in $K$, the bounds on $N$ and $N+M$ are linear in $C_M$. \qedhere
\end{proof}
 
\subsection{Optimal Order for Circle}\label{subsec:circle}
 
On $\mathbb S^1$ the functional $\Lam(Y,t)$ can be bounded in closed form, and it is also the value
of a linear program, so the circle allows both a proof and a numerical comparison. We first prove
the bound, then specialize to regular polygons, and finally compare both with computed values.
 
Write $\mathbb S^1=\mathbb R/2\pi\mathbb Z$. Then $\Lam(Y,t)$ is the value of the linear program
that maximizes $-\sum_jP(y_j)$ over trigonometric polynomials
$P=\sum_{m=1}^t(a_m\cos m\theta+b_m\sin m\theta)$ subject to $\frac1{2\pi}\int|P'|\le1$.
Rotating $Y$ changes neither $\Lam$ nor the spherical design property, so we may assume that no point of $Y$ equals $0$.  
For $Y=\{y_1,\dots,y_M\}$ let $\nu(s)=\#\{j:y_j\in[0,s]\}$ and $E(s)=\nu(s)-Ms/2\pi$, so that $E$
vanishes at both ends of $[0,2\pi)$. For a periodic $C^1$ function $P$ with $\int P=0$,
integration by parts gives
\begin{equation}\label{eq:circle-ibp}
  \sum_jP(y_j)=\int_0^{2\pi}P\dd \nu=\int_0^{2\pi}P\dd E=-\int_0^{2\pi}E(s)\,P'(s)\dd s,
\end{equation}
where the mean-zero condition removes the term $\frac{M}{2\pi}\int P$. Since $\int P'=0$, the
function $E$ may be replaced by $E-c$ for any constant $c$.
 
\begin{theorem}[Bound on the circle]\label{thm:circle}
Let $Y\subset\mathbb S^1$ be a spherical $t_1$-design with $M$ points. Then for every $t$,
including $t=\infty$,
\[
  \Lam(Y,t)\ \le\ \frac{\pi M}{t_1+1}.
\]
In particular, if $M\le C_Mt_1$ then $\Lam(Y,t)\le\pi C_M$ for all $t$.
\end{theorem}
The proof uses the Erd\H{o}s--Tur\'an inequality in Montgomery's form \cite[Ch.~1]{Mont94}: for $M$
points $y_j$ in $\mathbb R/2\pi\mathbb Z$ and every integer $H\ge1$, the discrepancy over arcs satisfies
\[
  \sup_I\Bigl|\#(Y\cap I)-\frac{M|I|}{2\pi}\Bigr|\le\frac{M}{H+1}+3\sum_{h\le H}\frac{|S_h|}{h},\qquad S_h=\sum_je^{ihy_j}.
\]
\begin{proof}[Proof of Theorem~\ref{thm:circle}]
Take $c$ to be the midpoint of $[\inf E,\sup E]$ in \eqref{eq:circle-ibp}. Then
\[
  \Bigl|\sum_jP(y_j)\Bigr|\le\tfrac12\operatorname{osc}(E)\int_0^{2\pi}|P'|\dd s
  =\pi\operatorname{osc}(E)\cdot\frac1{2\pi}\int|P'|.
\]
Because $E(s)-E(s')=\#(Y\cap(s',s])-M(s-s')/2\pi$, we have $\operatorname{osc}(E)\le\mathcal D(Y)$, where
$\mathcal D(Y)$ is the supremum over arcs $I$ of $|\#(Y\cap I)-M|I|/2\pi|$. As $Y$ is a spherical $t_1$-design,
$\sum_je^{ihy_j}=0$ for $1\le h\le t_1$. In the Erd\H{o}s--Tur\'an inequality with $H=t_1$ all
exponential-sum terms vanish, which leaves $\mathcal D(Y)\le M/(t_1+1)$. \qedhere
\end{proof}
 
\begin{remark}\label{rem:circle}
Theorem~\ref{thm:circle} does not use upper regularity of $Y$, which follows from the spherical design property and
$M\lesssim t_1$ by Proposition~\ref{prop:auto}. The theorem also removes the logarithm in the
case $d=1$ of Proposition~\ref{prop:dyadic}. Together with Proposition~\ref{prop:master}, it shows
that $N\ge\max\bigl(C_{1,d},\,\pi M/((t_1+1)c_d)\bigr)\,t$ new points suffice (with $d=1$), so for $M\le C_Mt_1$ the total is
$N+M\lesssim C_M\,t$, uniformly in $t_1<t$.
\end{remark}
 
The regular $M$-gon with $M=t_1+1$ is a tight spherical $t_1$-design with separation $2\pi/M$. For it
$\mathcal D(Y)=1$, and Theorem~\ref{thm:circle} gives $\Lam\le\pi$. This case also reduces to a problem in
one variable. Averaging $P$ over rotations by $2\pi/M$ leaves $\sum_jP(y_j)$ unchanged and does not
increase $\int|P'|$, so it suffices to consider $P(\theta)=Q(M\theta)/M$, where $Q$ is a mean-zero
trigonometric polynomial of degree at most $J=\lfloor t/M\rfloor$. After this rescaling,
\[
  \Lam(Y,t)=F(J),\qquad
  F(J)=\sup\Bigl\{Q(0):\ \deg Q\le J,\ \int Q\dd\mu=0,\ \frac1{2\pi}\int_0^{2\pi}|Q'|\le1\Bigr\},
\]
so $\Lam$ depends on $t$ and $t_1$ only through $J$. Three facts about $F$ follow.
\begin{enumerate}
  \item[(a)] $F(1)=\pi/2$. For $Q=r\cos(\theta-\varphi)$ one has $Q(0)\le r$ and
        $\frac1{2\pi}\int|Q'|=2r/\pi$.
  \item[(b)] $F(J)\le2H_J$, where $H_J=\sum_{k\le J}1/k$. This follows from $Q(0)=\sum_ma_m$ and
        $|a_m|\le2/m$, obtained by integrating by parts in $a_m=2\int Q\cos m\theta\dd\mu$.
  \item[(c)] $F(J)\le\pi$ for every $J$. Since $Q$ is periodic with mean zero,
        $\int_0^{2\pi}Q'(s)\,s\dd s=2\pi Q(0)$. Because $\int Q'=0$, this gives
        $Q(0)=\frac1{2\pi}\int_0^{2\pi}Q'(s)(s-\pi)\dd s$, and hence
        $|Q(0)|\le\pi\cdot\frac1{2\pi}\int|Q'|\le\pi$.
\end{enumerate}
Fact (c) is the case $\mathcal D(Y)=1$ of the computation in the proof of Theorem~\ref{thm:circle}. The bound
$2H_J$ in (b) is weaker than $\pi$ once $J\ge3$.
 
To compare these bounds with computed values, we solved the linear program with the gradient
constraint imposed at $G=8t$ equally spaced points, for two families with $t_1=8$. The first is the
regular $9$-gon. The second consists of two $9$-gons, one rotated by $0.3\cdot2\pi/9$. It is a
spherical $t_1$-design with $M=18$ that is upper regular with $K=2$, and Theorem~\ref{thm:circle} bounds its
$\Lam$ by $2\pi$.
 
\begin{table}[ht]
\centering
\begin{tabular}{rrcccc}
\hline
$t$ & $t/t_1$ & $J$ & regular $9$-gon & two rotated $9$-gons & $\min(\pi,2H_J)$\\
\hline
12 & 1.5 & 1 & 1.59 & 1.87 & 2.00\\
16 & 2   & 1 & 1.60 & 1.88 & 2.00\\
32 & 4   & 3 & 2.28 & 2.67 & 3.14\\
64 & 8   & 7 & 2.69 & 3.16 & 3.14\\
\hline
\end{tabular}
\caption{$t_1=8$, with $M=9$ for the polygon and $M=18$ for the pair, and $G=8t$. The last column is
the bound for the regular polygon from facts (b) and (c).}
\end{table}
 
The constraint is imposed only at grid points, so the computed values can
differ from $\Lam$. For $J=1$ the exact value is $\pi/2\approx1.571$, and the computed values
$1.59$ and $1.60$ are too large by $1\%$ to $2\%$. For $t_1=8$ and $t=32$, the grids $G=8t$ and
$G=24t$ give $2.281$ and $2.261$.
 
There are certain observations to be made here.  
\begin{enumerate}
  \item[(O1)] For the polygon, $\Lam$ depends on $t$ and $t_1$ only through $J$. At the same ratios $t/t_1$,
        the computed values for $t_1=4$ and $t_1=8$ agree within $3\%$, and those for $t_1=8$ and
        $t_1=12$ agree within $1\%$. The remaining differences are consistent with the differences in
        $J$ and with the discretization error.
        
  \item[(O2)] The computed values for the polygon rise from $1.59$ to $2.69$ as $t/t_1$ increases from
        $1.5$ to $8$ and stay below $\pi$, as fact (c) above requires. At $t/t_1=8$ they are about $86\%$ of
        the bound, so in this range the constant $\pi$ does not appear to be far from the true value
        for the polygon.
        
  \item[(O3)] Doubling the point set to two rotated polygons raises the computed values by a factor
        between $1.17$ and $1.18$. The bound of Theorem~\ref{thm:circle} allows a factor of $2$, and
        at $t/t_1=8$ the computed value is about $50\%$ of $2\pi$.
\end{enumerate}
 
It should be noted that only $t/t_1\le8$ was tested, and only for $d=1$, so the numerics say nothing directly about the exponent $t^{d-1}$ for $d\ge2$. A test in $d=2$ would require optimizing over polynomials of degree $t$ on $\mathbb S^2$ under an $L^1$-gradient constraint, and at feasible degrees the constants dominate. Such a test could also use the well-conditioned spherical designs of \cite{ACSW10} that Zheng and Zhuang propose as a numerical direction.
\section{Conclusion and Further Remarks}\label{sec:conclusion}

To conclude, the construction of \cite{ZZ} was kept unchanged, and the estimate of
$\Lam(Y,t)$ was improved. For a spherical $t_1$-design, $\Lam$ is controlled by an $L^1$ Jackson estimate and a sampling inequality applied on dyadic frequency blocks (Proposition~\ref{prop:dyadic}). For $d\ge2$ this gives order $t^{d-1}$, which cannot be improved for point sets of at most $c\,t^d$ points (Theorem~\ref{thm:sharpness}). For $d=1$ it gives the bound of Theorem~\ref{thm:circle}. The master estimate of Section~\ref{subsec:arbitrary} reduces the problem to the $L^1$ sampling constant of $Y$, and Proposition~\ref{prop:necessity} shows that some control of local multiplicity is needed.

Combined with Proposition~\ref{prop:auto}, the estimates above settle several questions that arise naturally from the extension problem. Items 1)--3)  depend on Proposition~\ref{prop:auto}, which uses the Bernstein inequality of Section~\ref{subsec:arbitrary}.
\begin{enumerate}
  \item[1)] \emph{Is bounded multiplicity at scale $1/t$ sufficient?} For spherical designs of optimal order the
        question does not arise. Applying Proposition~\ref{prop:necessity} to $Y$ itself shows that every
        spherical $t_1$-design with $M\lesssim t_1^d$ has bounded multiplicity at scale $c_*/t_1$, which is a
        stronger property and gives upper regularity without further hypotheses.
        
  \item[2)] \emph{Do arbitrary optimal-order spherical designs extend with $O(t^d)$ points?} Yes, with a constant
        depending on $C_M$. By Theorem~\ref{thm:arbitrary}, every spherical $t_1$-design with $M\le C_Mt_1^d$ extends
        to a spherical $t$-design with $N+M\le C(d)\,C_M\,t^d$, uniformly in $t_1<t$ and without a separation
        hypothesis. The linear dependence on $C_M$ cannot be removed in general, since a spherical design made of
        $m$ copies of a smaller spherical design needs at least $c\,m\,t^d$ points by
        Proposition~\ref{prop:necessity}. For $M\lesssim t_1^d$ this improves the bound $t^{2d}$ of
        Corollary~\ref{cor:arbitrary} and the bound $t^{2d+1}$ of \cite{ZZ}.
        
  \item[3)] \emph{Is there an optimal-order spherical design that forces more than order $t^d$ new points?} No, with
        constants depending on $C_M$. By the previous item, the classical bound $N+M\ge ct^d$ is the only
        constraint on the order of the total for such designs.
        
  \item[4)] \emph{What happens on the circle?} Theorem~\ref{thm:circle} gives $\Lam\lesssim1$ for every
        spherical $t_1$-design with $M\lesssim t_1$, without a logarithm.
\end{enumerate}
Item 2) gives an exact completion of order $t^d$ for every optimal-order spherical $t_1$-design, uniformly in $t_1<t$. 
The preprint \cite{NestedQMC} by Wu and Zhuang studies nested equal-weight QMC designs on spheres, that is, point sets with the optimal Sobolev rate that need not be polynomially exact, and presents them as a counterpart to the nested designs of \cite{ZZ}. Among other results, it shows that every prescribed finite point set admits a completion to a QMC design of optimal rate, using a step that compensates for the low-frequency discrepancy of the inherited nodes. These results concern integration error rather than exact polynomial integration, so they are distinct from Theorem~\ref{thm:arbitrary}, which gives an exact completion for every starting spherical design of optimal order. Our method belongs to the line of work that began with \cite{BRV13}: degree theory produces a zero of a vector field on polynomial space, and all quantitative content is in estimates for linear functionals on $\partial\Omega$. In that sense, our contribution is an estimate of one such functional that uses the design property of the given points together with their local regularity.

There are certain questions that future work could take up. \begin{enumerate}
  \item[F1)] The conjecture of \cite[Section~5]{ZZ} concerns the constants in the
        orders $t_1^d$ and $t^d$. Our constants depend on $d$ and $C_M$ and, through the extension step,
        on $C_{1,d}$ and $c_d$, and we do not compare them with the constants of the starting design.
  \item[F2)] Theorem~\ref{thm:chain} has constants of size
        $\gamma_d^{\,L}$. Proposition~\ref{prop:auto} allows the propagation lemma to be replaced by a size
        bound $|Y_k|\le C_kt_{k+1}^d$, but the multiplicity constant is proportional to $C_k$, and the
        number of new points is proportional to that constant, so the geometric growth remains. A bound
        independent of $L$ may require the new points to be separated, as in \cite{BRV15}.
  \item[F3)] The values of $C_{1,d}$ and $c_d$ come from the
        construction in \cite[Appendix~B]{ZZ} and are probably not optimal. Since $\Lam\asymp t^{d-1}$
        for upper regular starting designs, the gap between the number of new points certified by the Brouwer
        scheme and the classical lower bound depends only on these constants.
  \item[F4)]  For optimal-order spherical designs the order of $N+M$ is now
        determined up to constants, so a lower bound for the nested problem beyond $N+M\ge ct^d$
        could only come from spherical designs with $M\gg t_1^d$. Proposition~\ref{prop:necessity} already covers
        clusters, and it is not known whether clusters are the only obstruction for such designs.
  \item[F5)]  The test in Section~\ref{subsec:circle} covers only $d=1$
        and $t/t_1\le8$. A study for $d=2$, for instance with the designs of \cite{ACSW10}, would show
        how large the constants are in practice.
\end{enumerate}


\end{document}